\documentclass[11pt]{amsart}   %%%%%%%Please do not change

       \theoremstyle{plain}
       \newtheorem{theorem}{Theorem}
       \newtheorem{lemma}[theorem]{Lemma}

       \newtheorem{definition-theorem}[theorem]{Definition}

      \newtheorem{definition}{Definition}[theorem]
      \newtheorem{remark}{Remark}

       \newcommand{\C}{{\mathbb C}}
       
       \newcommand{\Z}{{\mathbb Z}}

       \newcommand {\stre} {\operatorname{tr_s}\exp}

       \makeatletter
       \def\@setcopyright{}
       \def\serieslogo@{}
       \makeatother

\begin{document}
 \title[Residue Chern Character]{A Note On the Residue Chern Character.}

    \author{Dmitry Gerenrot}
    \address{ Dmitry Gerenrot, Georgia Institute of Technology,
    School of Mathematics, 686 Cherry Street, Atlanta,
    GA, 30332-0160, U.S.A.}
    \email{gerenrot@math.gatech.edu}

%    \curraddr{School of Mathematics, Georgia Institute of Technology,
 %     , State College PA 16802}

\thanks{ Manuscript received (to be filled by editor); revised
    (to be filled by editor).
       }

    % abstract (optional)
\begin{abstract}
The aim of this note is to improve upon our earlier result
\cite{Ger} which translates Weyl's (curvature) formulation of
Chern character of a smooth vector bundle  (\cite{MS,Q}) 
into the language of residues. Given super-vector bundle $E$ over a compact manifold
$M$, one pulls it back to $T^*M.$ The dualized Chern character is
the functional $\eta\mapsto\int_{T^*M}ch(\pi^*E)\pi^*\eta$ which
is expressed in \cite{Ger} as a sum of the residues of the
form $Res_{s=z}\Gamma(s) \int_{T^*M}\nabla^{-2s}_L\pi^*\eta$ over
all $z\in {\C}.$ Here, $\nabla^2_L$ is constructed from a
curvature of the pullback bundle $\pi^*E$ to $T^*M$ and
an odd endomorphism $L$ of $\pi^*E$ using Quillen's formalism
\cite{Q}. The result in our previous paper is stated using a limit
as ${R\to 0}$  of the integrals over $T^*M$ minus the $R$-tubular
neighborhood of the zero section. The present note 
provides for a stronger, more effective  formulation.
\end{abstract}

    % AMS keywords (used in AMS journals)
    \keywords{Noncommutative Geometry Chern Character Connes Moscovici
            Quillen Superconnection Residue Cocycle}

    % today's date, or fill in whatever date you prefer
    \date{\today}

    \maketitle
%\section{gftg}

Our result in \cite{Ger}
resembles the Local Index Formula (\cite{CM, Hig}) in a
differential-geometric setting. Let $M$ be a smooth $n$-manifold
with no boundary and let $E$ be a smooth ${\Z}_2$-graded vector
bundle. Let $\pi:T^*M\to M$ be the cotangent bundle and let $L$ be
an odd skew-adjoint endomorphism of $\pi^*E$ invertible everywhere
but the zero section of $T^*M$. Finally, we assume that the
coefficients of $L$ are first-order homogeneous polynomials of the
"vertical" (fiberwise) coordinates of $T^*M.$ For example, $L$
could be a symbol of an odd self-adjoint elliptic operator on $E$.

Now, let $\nabla$ be a connection on $\pi^*E$ which is a pullback
of some connection on $E$. Suppose that both connections respect
the grading of $E$ and $\pi^*E$. According to Quillen \cite{Q},
the Chern character corresponding to $L$ may be written as
$\operatorname{exp}(\nabla+L)^2$. Here, $\nabla+L$ is a differential operator on
the sections of $\Lambda^* T^*M\otimes\pi^* E$ which shares many
properties of an ordinary connection, such as the ${\Z}_2$-graded
Leibniz rule and the fact that $(\nabla+L)^2$ is an endomorphism
of $\Lambda^* T^*M\otimes\pi^* E$, i.e. a $0$-th order
differential operator, which is still called the {\it curvature}.
In addition, the graded trace (the {\it supertrace})
${\operatorname tr_s}(\nabla+L)^{2k}$ for any $k$ is a
differential form whose cohomology class is independent on the
concrete choice of $\nabla$, a topological characteristic of $E$.
We denote $\nabla+L$ by $\nabla_L$ and call it a {\it
superconnection}.

The big advantage is that if $L$ is as above, then the Chern
character form ${\stre}\nabla^2_L$ decays exponentially fast along
the fibers of $T^*M$ so that the dual Chern character current can
be defined on $\Omega^*M$:
$$\eta\mapsto\int_{T^*M}\pi^*(\eta){\stre}\nabla^2_L.$$
WARNING: for notational convenience, we omit the supertrace from
our formulas, though it is tacitly assumed everywhere.

In \cite{Ger}, we have proved the following theorem about this current:

\begin{theorem}\label{theorem:oldthm}
Let $Y_R$ be the open $R$-tubular neighbourhood of the zero 
section in $T^*M$ and let $X_R$ be its complement.
Under the hypotheses outlined above, for any $\eta\in\Omega^*(M)$
\begin{multline}\label{equation:oldthmeq}
\int_{T^*M}\!\!\!\!\operatorname{tr_s}\pi^*(\eta)\exp{\nabla_L^{2}} =
\lim_{R\to 0}\sum_{z\in{\C}}Res|_z\Gamma(z)
\int_{X_R}\operatorname{tr_s}\pi^*(\eta)({-\nabla_L^2})^{-z},
\end{multline}
where the right-hand side integral is understood to be the 
meromorphic extension from the region $Re(z)\gg 0,$ on which
it converges.
Further, all but finitely many residues on the right-hand side
vanish as $R\to 0$.
\end{theorem}

The proof of the equality \ref{equation:oldthmeq} is based on the Mellin
Transform. Section 6 of \cite{Ger}, contains the
argument and a concise outline preceding it.  In the present note,
we concentrate on the right-hand side of the equality. 
We state the 
following slightly stronger result.

\begin{theorem}\label{theorem:newthm}
Under the hypotheses outlined above, for any
$\eta\in\Omega^\kappa(M)$ and any positive $R,$
\begin{align}\label{equation:newthmeq}
\int_{T^*M}\!\!\operatorname{tr_s}\pi^*(\eta)\exp{\nabla_L^{2}}\!\! &=
Res_{z={\frac \kappa
2}-n}\Gamma(z)\int_{X_R}\operatorname{tr_s}\pi^*(\eta)\big[({-\nabla_L^2})^{-z}\big]_{2n-\kappa},
\end{align}
where the right-hand side integral is understood to be the 
meromorphic extension from the region $Re(z)\gg 0,$ on which
it converges. In particular, it does not depend on $R.$
Further, both sides vanish if $\kappa$ is odd.
\end{theorem}

Here, by $[\omega]_\kappa$ we denote the $\kappa$-degree
part of the mixed differential form $\omega.$ 
We proceed to:
\begin{itemize}
\item[1)] Review the notion of complex powers $(-\nabla_L^2)^{-z}$
via holomorphic functional calculus. \item[2)] Review the geometric series 
expansion of $(-\nabla_L^2)^{-z}$ used in \cite{Ger}. \item[3)] Prove theorem
\ref{theorem:newthm} based on that expansion and on theorem \ref{theorem:oldthm}.
\end{itemize}

1) Complex powers of the curvature $\nabla_L^2$ are defined via
the following integrals:
$$(-\nabla_L^2)^{-z}={\frac 1 {2\pi i}}\int_\gamma \lambda^{-z}
(\lambda+\nabla^2)^{-1}d\lambda,$$ where $\gamma$ is a counter-clockwise
oriented contour which surrounds the pointwise spectrum of
$\nabla_L^2$. We prove in \cite{Ger}, section 5, that 
$\gamma$, in fact,  may be taken as a vertical which is oriented 
downward and separates $sp(\nabla_L^2)$ from the imaginary axis. 
Such $\gamma$ exists as long as the underlying point of $T^*M$ 
does not lie in the zero section.  In \cite{Ger}, we  have also dealt
with the fact that $\gamma$ depends on that point in the first place.
(Briefly, we have shown that if we integrate over $X_R$, then
$\gamma$ can be chosen uniformly. But then we take the limit 
as $R\to 0$.) We have also
shown that $\int_{X_R}\pi^*(\eta)(-\nabla_L^2)^{-z}$ converges
for $Re(z)\gg 0$ and has a meromorphic extension to all of
${\C}$ with at most simple poles.

2)In order to see why the meromorphic extension exists, we write
out the expression for $(-\nabla_L^2)^{-z}$ via geometric series. First,
\begin{align*}
\nabla_L^2&=(\nabla+L)^2\\
&=\nabla^2+[\nabla,L]+L^2\\
&=\nabla^2+dL+[\theta,L]+L^2,
\end{align*}
where $\nabla=d+\theta$ and $\theta$ is a locally defined odd
endomorphism of $\pi^*E$. Such $\theta$ exists, for any connection
or superconnection and is comprised of Christoffel symbols.
"Locally" means that we restrict our attention to a coordinate
chart $V$ on $T^*M$ over which $\pi^*E$ is trivial and which is 
itself a local trivialization of $T^*M$ over $M$. Thus, there are two sets
of coordinates. "Horizontal", i.e. coordinates of $M$ $x^1\ldots
x^n$; and "vertical" coordinates of a fiber. We use polar
coordinates here, $\rho$ being the radial one and
$\Xi^1,\ldots\,\Xi^{n-1}$ being the coordinates of a unit sphere.
We then write:

\begin{align}
\label{equation:weget1}{2\pi i}(-\nabla_L^2)^{-z}\!\!&=\!\!\int_\gamma\!\!\scriptstyle
\lambda^{-z}\big(\lambda+L^2 +d_x L+d_\Xi L
+[\theta,L]+d_\rho L+\nabla^2\big)^{-1}\!\!d\lambda\\ \nonumber
\!\!&=\!\!\int_\gamma\scriptstyle \!\!\lambda^{-z}(\lambda+L^2)^{-1}
\sum_{k=0}^{2n}\big(-(\lambda+L^2)^{-1} (d_x L+d_\Xi L
+[\theta,L]+d_\rho L+\nabla^2)\big)^k \!\!d\lambda.
\end{align}

Observe that $\nabla^2+d_\rho L$ is not a multiple of $\rho$, while
from $d_x L+d_\Xi L +[\theta,L]$ one power of $\rho$ may be pulled
out. We then expand the $k$-th term of the series as a
non-commutative polynomial in $\nabla^2$, $d_\rho L$, $d_\Xi L$,
and $d_x L+[\theta,L]$ and separate the powers of $\rho$.

For illustration, we treat just one typical "interesting" term. It
must contain one copy of $d_\rho L$ and $n-1$ copies of $d_\Xi L$
in order to produce a volume form on $T^*M$. (Observe that
$\theta$ and $\nabla^2$ cannot contain any vertical differentials,
since $\nabla$ has been pulled back from $E$.) One such term has
the form:
\begin{align}
\label{equation:weget2}\nonumber\int_\gamma\lambda^{-z}&(\lambda+L^2)^{-1}
\!\big[(\lambda+L^2)^{-1} (d_x L+\![\theta,L])\big]^l\,\times
 \\ \nonumber&\qquad
\big[(\lambda+L^2)^{-1} d_\rho L\big] \big[(\lambda+L^2)^{-1}
d_\Xi L\big]^{n-1}
 \big[(\lambda+L^2)^{-1}\nabla^2)\big]^{k-n-l} \!d\lambda
\\ \nonumber&
=\rho^{-2(z+k)+n+l-1}
 \int_\gamma\sigma^{-z}(\sigma+{\L}^2)^{-1}
\big[(\sigma+{\L}^2)^{-1} (d_x {\L}+[\theta,{\L}])\big]^l\,\times
\\ &\qquad
\big[(\sigma+{\L}^2)^{-1} d_\rho L\big] \big[(\sigma+L^2)^{-1}
d_\Xi {\L}\big]^{n-1}
 \big[(\sigma+{\L}^2)^{-1}\nabla^2)\big]^{k-n-l} d\sigma,
\end{align}
where ${\L}={\frac L \rho}$ and $\sigma={\frac \lambda {\rho^2}}$.
Strictly speaking $\gamma$ should be replaced by ${\frac 1
{\rho^2}} \gamma$ but we leave that detail to \cite{Ger}.

3) In order to prove theorem \ref{theorem:newthm}, we need to compute
$\int_{X_R\bigcap V}\pi^*(\eta)(-\nabla^2_L)^{-z}$.  Rather than
using $(-\nabla_L^2)^{-z}$, we perform the computation just for
the sample term shown in (\ref{equation:weget2}). Assuming $Re(z)\gg 0$, we
get:
\begin{align}\label{equation:weget3}
\nonumber
\int_{X_R\bigcap
V}&\!\!\!\pi^*(\eta)\!\int_\gamma\!\lambda^{-z}(\lambda+L^2)^{-1}
\big[(\lambda+L^2)^{-1} (d_x L+[\theta,L])\big]^l\,\times
\\ &\nonumber \,\,
\big[(\lambda+L^2)^{-1} d_\rho L\big] \big[(\lambda+L^2)^{-1}
d_\Xi L\big]^{n-1}
 \big[(\lambda+L^2)^{-1}\nabla^2)\big]^{k-n-l} d\lambda
\\ & \nonumber\\ \nonumber&
=\int_R^\infty \pi^*(\eta)\rho^{-2(z+k)+n+l-1 }d\rho \int_{\pi(V)}dx \int_{S^{n-1}} d\Xi
 \,\times
\\ \nonumber &
\qquad\qquad\int_\gamma\sigma^{-z}(\sigma+{\L}^2)^{-1}
\big[(\sigma+{\L}^2)^{-1} (d_x {\L}+[\theta,{\L}])\big]^l\times
\\ \nonumber &\qquad\qquad\qquad
\big[(\sigma+{\L}^2)^{-1} d_\rho L\big] \big[(\sigma+L^2)^{-1}
d_\Xi {\L}\big]^{n-1}\times
\\ \nonumber &\qquad\qquad\qquad
 \big[(\sigma+{\L}^2)^{-1}\nabla^2)\big]^{k-n-l} d\sigma
\\&\nonumber\\&
= {\frac {R^{-2(z+k)+n+l}}{2(z+k)-n-l}}\phi_V(z),
 \end{align}
 where $\phi_V$ is defined by the above equation. It is entire in
 $z$. (It is true after $\sigma$ is integrated out and
 then still true after one integrates out $\Xi$ and $x$ over the compact
 manifold $S^*M$\cite{Ger}). Also, $\phi_V$ is independent on $R$. 
Thus, the integral in \ref{equation:weget2} has a meromorphic extension to ${\C}$
with at most simple poles.
Compactness of $M$
 and the fact that there are only finitely many
 terms such as the one above, imply that
 $\int_{X_R}\pi^*(\eta)(-\nabla_L^2)^{-z}$ also has a meromorphic
 extension to ${\C}$.

We proceed to prove the theorem \ref{theorem:newthm}. For the case
of even $\kappa=deg(\eta),$ we consider the residues of:
$$\Gamma(z){\frac {R^{-2(z+k)+n+l}}{2(z+k)-n-l}}\phi_V(z).$$
We have the following lemma:

\begin{lemma}\label{theorem:l1}
$\phi_V(-m)=0 $ for all $m < {\frac {3n-\kappa}2}.$
\end{lemma}
{\bf Proof:} Observe that for $z=0,-1,\ldots -(k+1)$ $\phi_V(z)=0$,
since then $(-\nabla_L^2)^{-z}$ is just a positive integer power.
To see this, suppose for a moment that $L^2$ is a scalar and thus:
$$(-\nabla_L^2)^{-z}=\sum_{k=0}^\infty \big(^{-z}_k\big)(L^2)^{-z-k}
                                       ([L,\nabla]+\nabla^2)^k.$$
If $z=-m$, the series terminates for $k>m$. For general $L$,  the phenomenon 
is similar. Next,
we count the differential form degrees in (\ref{equation:weget3}).
Unless
$$\kappa+2(k-n-l)+l+n=2n=dim(T^*M),$$  $\phi_V$ vanishes identically. Thus,
$\kappa+2k-l=3n$, so that $${\frac {3n-\kappa}2}\le k.$$\qed

Now, one can easily compute the location of the residue due to
${\frac {R^{-2(z+k)+n+l}}{2(z+k)-n-l}}$ in terms of $\kappa$:
$${\frac {n+l} 2}-k={\frac \kappa 2}- n.$$
In summary, all the residues coming from $\Gamma(z)$ between $0$
and ${\frac {\kappa-3n}2}$ are killed by the zeroes of $\phi_V$
and the only residue in that range is supplied by $${\frac
{R^{-2(z+k)+n+l}}{2(z+k)-n-l}}= {\frac
{R^{-2(z+n)-\kappa}}{2(z+n)-\kappa}}.$$

Finally, in \cite{Ger} we have shown that 
\begin{multline}\label{equation:finally}
\int_{X_R}\!\!\!\!\operatorname{tr_s}\pi^*(\eta)\exp{\nabla_L^{2}} =
\sum_{z\in{\C}}Res|_z\Gamma(z)
\int_{X_R}\operatorname{tr_s}\pi^*(\eta)({-\nabla_L^2})^{-z},
\end{multline}
The argument uses the fact that $\int_0^\infty e^{-\lambda t}t^{z-1}dt=
\lambda ^{-z}\Gamma(z)$ and Mellin transforms to translate the 
exponential expression on the left of (\ref{equation:finally}) 
into the one on the right.
Taking limits as $R\to 0$, we see that the residues to the left of
${\frac {\kappa-3n}2}$ are multiples of a positive power of $R$
and thus vanish, whereas the one at ${\frac \kappa 2}-n$ is
independent on $R$. This proves the theorem 
(\ref{equation:newthmeq}) in the case when
$\kappa$ is even. 

If $\kappa$ is odd, then due to $\kappa+2k-l=3n$,  $l+n$ must be
odd. That is the total power of $d_\rho L$, $d_\Xi L$ and $d_x
L+[\theta,L]$. Now, since all the connections preserve the
${\Z}_2$-gradings, the local expression for $\theta$ must be a
block-diagonal matrix of 1-forms. The same is true about
$\nabla^2$. However, $L$ is block-off-diagonal, since it is an odd
endomorphism. Thus, if $n+l$ is odd, the supertrace of the
corresponding term vanishes, which proves the theorem \ref{theorem:newthm} for
odd $\kappa$.

\begin{remark}
{\rm It is clear that 
the right side of theorem \ref{theorem:newthm} may be written as 
an integral over the unit sphere bundle. It would be curious to 
obtain this through Wodzicki Residue.}
\end{remark}

\end{document}